\documentclass[11pt]{amsart}

\usepackage{a4wide,amssymb,amsmath}

\usepackage{upref}
\usepackage{mathtools}
\usepackage{verbatim}
\usepackage{amscd}

\newcommand{\NN}{\mathbb{N}}
\newcommand{\RR}{\mathbb{R}}
\newcommand{\QQ}{\mathbb{Q}}
\newcommand{\ZZ}{\mathbb{Z}}
\newcommand{\id}{\mathrm{id}}
\newcommand{\abs}[1]{\left|#1\right|}

\newcommand{\fix}{\mathrm{fix}}
\newcommand{\orb}{\mathrm{orb}}
\newcommand{\relReal}{\mathcal{ER}_{\mathrm{rel}}}

\newtheorem{theorem}{Theorem}[section]

\newtheorem{prop}[theorem]{Proposition}
\newtheorem{corollary}[theorem]{Corollary}
\theoremstyle{definition}
\newtheorem{rem}{Remark}
\newtheorem{example}{Example}

\begin{document}
\title[Relative realizability]
{Realizability of integer sequences\\ \vspace{1mm} as differences of
fixed point count sequences}
\author{Natascha Neum\"arker}
\address{Fakult\"at f\"ur Mathematik, Universit\"at Bielefeld, 
Postfach 100131, 33501 Bielefeld, Germany}
\email{naneumae@math.uni-bielefeld.de}

\begin{abstract}
A sequence of non-negative integers is \textit{exactly realizable}
as the fixed point counts sequence of a dynamical system if and only 
if it gives rise to a sequence of non-negative orbit counts. 
This provides a simple realizability criterion based on 
the transformation between fixed point and orbit counts.
Here, we extend the concept of exact realizability to realizability 
of integer sequences as differences of the two fixed point counts 
sequences originating from a dynamical system and a topological factor.
A  criterion analogous to the one for exact realizability is given 
and the structure of the resulting set of integer sequences is outlined.  
\end{abstract}
\maketitle


\section{Introduction and general setting}
A (topological) dynamical system $(X,T)$ is given by a 
topological space $X$ and a continuous map $T: X\to X$.
Associated to each dynamical system are two sequences
$a=(a_n)_{n\geq 1}$ and $c=(c_n)_{n\geq 1}$,
defined by
\begin{equation}
    a_n = \left|\{ x\in X\mid T^n (x) = x\}\right|\ \text{and}\ 
    c_n = \left|\{\mathcal{O}\mid \mathcal{O}\ \text{is a periodic
		orbit of length}\ n\}\right|.
\end{equation}
They count the fixed points and periodic orbits,
respectively.
In interesting classes of dynamical systems, these numbers are
finite for all $n$, which we assume from now on.
Thus, $a$ and $c$ are sequences of non-negative
integers.

Clearly, $a$ and $c$ are related by M\"obius inversion,
\begin{equation}\label{E:fixOrb}
     \fix(c)_n\coloneqq a_n = \sum_{d|n} d\, c^{}_d \quad
			      \text{and}\quad 
     \orb(a)_n\coloneqq c_n = \frac{1}{n}
			\sum_{d|n}\mu\left(\frac{n}{d}\right) a^{}_d.
\end{equation}
Details on these transformations, considered as linear operators
on the space of arithmetic functions, 
can be found in \cite{fixOrb} and references given there.
In what follows, we sometimes make use of the operator notation
$\orb(a)$ and
$\fix(c)$, referring to the transformations of the sequences $a$ and
$c$  according to \eqref{E:fixOrb}.
Since $\fix$ and $\orb$ are inverses of each other, in the sense that
$\fix\circ\orb = \orb\circ\fix=\id$, it is possible
to implicitly define a sequence $f$ by setting
$\orb(f) = g$ or $\fix(f)=g$ for a given sequence $g$.
Note that $\fix$ and $\orb$ are well-defined for arbitrary sequences
of complex numbers.
If, however, $f$ is a sequence of integers, then so is $\fix(f)$,
whereas the converse need not be true. 

According to Puri's terminology \cite{purisDiss},
a \textit{system} $(X,T)$ comprises an arbitrary set $X$ and a map
$T:X\to X$.
An integer sequence $(f_n)_{n\geq 1}$ is called \textit{exactly
realizable} if there is a system $(X,T)$ whose fixed point counts are
given by $(f_n)_{n\geq 1}$. Membership in the set $\mathcal{ER}$ of
exactly realizable sequences is characterized by the Basic Lemma
\cite[Thm.2.2]{purisDiss}:
\begin{theorem}\label{basicLemma}
    A  sequence of non-negative integers $(f_n)_{n\geq 1}$ is exactly
    realizable if and only if, for all $n\geq 1$, the sum
    $\sum_{d|n}\mu\left(\frac{n}{d}\right) f^{}_d$ yields a
    non-negative integer divisible by $n$, that is, if
    and only if $\orb(f)$
    is a sequence of non-negative integers.\hfill 
\end{theorem}
In other words, the only restriction for an integer sequence to be
exactly realizable is that it gives rise to an orbit counts sequence
of non-negative integers. By compactification of $\NN$ with respect 
to the discrete topology and the definition of an appropriate
permutation of the resulting set, Puri shows that the realizing
system can be chosen to be a homeomorphism on a compact space.
Windsor \cite{windsor} even gives a construction of a
smooth system realizing an arbitrary sequence from $\mathcal{ER}$.
Moss \cite{mossDiss} systematically investigates the realization of integer sequences
by algebraic dynamical systems. 
A collection of many exactly realizable sequences from the
Online Encyclopedia of Integer Sequences (OEIS) \cite{OEIS}
is listed by Puri and Ward \cite{PuriWard}.

A \textit{(topological) factor} of a dynamical system $(X,T)$
is a dynamical system $(Y,S)$ for which there is a 
continuous surjection
$\phi: X\to Y$ that makes the diagram
\begin{equation}\label{CD}
    \begin{CD}
       X            @>T>>   X\\
       @V{\phi}VV   @VV{\phi} V\\
       Y            @>S>>   Y\\
    \end{CD}
\end{equation}
commutative, i.e.,
$\phi(T(x)) = S(\phi(x))$ for all $x\in X$
and, by induction, $\phi(T^n (x)) = S^n (\phi(x))$.

Being a factor of the dynamical system $(X,T)$ is a much weaker
condition than topological conjugacy, for which $\phi$ is 
required to be a homeomorphism.
While the orbit statistics of both systems 
coincide in the latter case, a factor can have a completely different
orbit structure.
Going down to the factor, an arbitrary number of periodic points can 
be `lost' since, for some arbitrary map $U$, the map $U\times S$ 
always gives rise to the factor $S$.
An example of the number of periodic points in the factor system
exceeding the number of periodic points in the original system
is given by the maps dual to $x\mapsto 2x$ on $\QQ$ and $\ZZ$,
respectively, see example \ref{SIntegers}.

Since for many dynamical systems it is known that they are a
factor of some well-understood dynamical system or, vice versa,
that some well-studied system is a factor (cf.\ \cite{torusParam}),
the question is raised,
in what way their orbit statistics can be related.
A natural approach to this question is the classification of
integer sequences as \textit{relative} fixed point counts, that is,
as the difference of the fixed point counts sequences asscociated to
$(X,T)$ and $(Y,S)$.
It is clear, however, that this is only a very coarse way of
relating the dynamics of the two systems. 
Example~\ref{TMPD} later on gives an example of a dynamical system
and a factor whose fixed point counts coincide, the factor map being
stricly $2:1$, however.

The purpose of this short note is to
derive  a realizability criterion for integer sequences
as such difference sequences,
analogous to the one given in Theorem \ref{basicLemma},
and to consider a few consequences.

\section{Relative realizability}
An integer sequence $h$ is called \emph{relatively realizable} if there
is a dynamical system $(X,T)$  and a factor $(Y,S)$ 
with fixed point counts sequences $f$ and $g$, respectively, such that
$h = f - g$.
We denote the set of relatively realizable sequences by
$\relReal$. 

By linearity of the mapping $\orb$, the orbit counts sequences
corresponding to $f$ and $g$ 
satisfy the same relation: 
$\orb(h) = \orb(f) - \orb(g)$.
Since the realization is constructed orbit-wise, it is thus more
convenient to consider orbit counts. 

If the periodic orbits of $(X,T)$ and $(Y,S)$ partition the
respective space, each surjective map on the set of equivalence classes
defined by the periodic orbits can be extended to a
surjection $\phi:X\to Y$. 
More precisely, if a $T$-orbit 
$\mathcal{O} = \{ x, Tx ,\ldots , T^{n-1} x\}$
is assigned an $S$-orbit $\{y, S y,\ldots, S^{d-1} y\}$, the
corresponding map $\phi:X\to Y$ is defined via
\begin{equation}\label{elementwise}
    \phi(T^k(x)) = S^{r_k}(y),\quad  
	r_k\leq d-1,\quad r_k  \equiv k\ \mod d
\end{equation}
for the elements of $\mathcal{O}$.

\begin{prop}\label{powerSeriesCrit}
    A sequence of integers $h$ is relatively realizable if and only if
    the formal power series 
    $H(x) \coloneqq \sum_{n=1}^{\infty} h_n x^n\in\ZZ[[x]]$ 
    associated to $\orb(h)$ admits a decomposition
    \begin{equation*}
	H(x) = \sum_{n=1}^{\infty} b_n x^n + \sum_{n>d|n} 
		    a^{}_{d,n} (-x^d + x^n),
    \end{equation*}
    with $a^{}_{d,n}, b_n\in\NN_0$. 
\end{prop}
\begin{proof}
The grouping of the terms in $H(x)$ encodes the definition of an
appropriate surjection $\phi$.
We first note that the $n$-th term has coefficient
\begin{equation}\label{nthCoeff}
    b^{}_n + \sum_{n > d | n} a^{}_{d,n} - \sum_{k\geq 1} a^{}_{n,kn}.
\end{equation}
Let $(X,T,Y,S,\phi)$ relatively realize $h = f - g$ and set
$\orb(f) = (\nu_1,\nu_2,\ldots)$,
$\orb(g)=(\gamma_1,\gamma_2,\ldots)$.
For each $\nu_k \geq 0$, let $C_1,\ldots, C_{\nu^{}_k}$ denote the
$T$-orbits of length $k$.
For $d|k$ with $d<k$ set
\begin{equation*}
    a_{d,k} = \left|\left\{ D\in\{\phi(C_{1}),\ldots,
	    \phi(C_{\nu_k})\}: \right.\right.
	    \left.\left. \abs{D}=d, D\ \text{has no preimage orbit of
	    length}\ < k\right\}\right|.
\end{equation*}
In particular, $a_{d,k} = 0$ for divisors $d$ of $k$ that do not show
up as the length of an
image cycle and coinciding image cycles $\phi(C_i)=\phi(C_j)$ for
$i\not= j$ are counted only once. 
Then, the number of $\ell$ with
$a_{d,\ell d} \not= 0$ is finite, and 
$\sum_{\ell\geq 1} a_{d,\ell d} = \gamma_{d}$.
Define
$b_k \coloneqq \nu_k - \sum_{n > d|n}^{} a_{d,k}$ if
$\nu_k > 0$,  giving $b_k = 0$ otherwise.
Thus, the $a_{d,k}, b_k$ define a formal power series of the
structure
indicated above whose $n$-th coefficient is, according to
\eqref{nthCoeff},
\begin{equation*}
    b^{}_n + \sum_{n > d | n} a^{}_{d,n} - \sum_{k\geq 1} a^{}_{n,kn}
    = \nu_n - \gamma_n = \orb(h)_n.
\end{equation*}
Consequently, this is an appropriate decomposition of the power
series of $\orb(h)$.

For the inverse direction, we define
\begin{equation*}
    \nu_n \coloneqq b_n + 1 +  \sum_{n>d|n} a_{d,n}\ \text{and}\
    \gamma_n \coloneqq \sum_{k\geq 1} a_{n,kn} + 1
\end{equation*}
and realize the sequences $\orb(f) \coloneqq (\nu_1,\nu_2,\ldots)$ and
$\orb(g) \coloneqq (\gamma_1,\gamma_2,\ldots)$ on the one-point
compactification $(\NN_*,\tau_*)$ of $\NN$.
If `the point at infinity' $\infty$, added for the purpose of
compactification, is a fixed point of the permutation 
$\sigma$ on $\NN_* = \NN\cup\{\infty\}$, $\sigma$ is continuous in the
obtained topology,
see \cite{purisDiss} and the reference given therein.
In each of the resulting systems $(\NN_*,T)$ and $(\NN_*,S)$, there is 
at least one cycle of length $n$ for all $n\in\NN$.
By construction, $(\NN_*,S)$ and $(\NN_*,T)$ have $a_{n,kn}$
corresponding cycles of length $n$ and
$kn$, respectively. $\phi$ can be defined on these as described by
\eqref{elementwise}.
For the $b_n$ further  $T$-cycles of length $n$, one of the 
$n$-$S$-cycles can be chosen as the image cycle. 
If  the $1$-orbit $\{\infty\}$
is mapped to its analogue in the factor system, it is straightforward
to check that the map $\phi$ is continuous
with respect to the considered topologies. It is surjective by
construction, thus turning $(\NN_*,S)$ into a factor of $(\NN_*,T)$ 
and therefore
yielding a relative realization of the sequence $h$.
\end{proof}
As a consequence, we obtain that
each integer sequence is the difference orbit counts sequence of a
dynamical system and a factor.
\begin{theorem}\label{BLanalogue}
    A sequence $h$ of integers is relatively realizable if and
    only if $\orb(h)$ is a sequence of integers.
\end{theorem}
\begin{proof}
Since it is clear that the condition is necessary, it suffices to
give a decomposition of an arbitrary element 
$H(x)=\sum_{k=1}^{\infty} \eta^{}_k x^k $
from $\ZZ[[x]]$
as in Proposition~\ref{powerSeriesCrit}.
A possible approach for doing so is to 
select a divisor $d|n$ and a multiple $\ell n$ of $n$
and to decompose $\eta^{}_n$ 
into a sum of the shape
$\eta_n = b_n + a^{}_{d,n} - a^{}_{n,\ell n}$
with non-negative integers $b_n, a^{}_{d,n}$
and  $a^{}_{n,\ell n}$. 

Consider the case $n$ odd first.
If $\eta^{}_n \geq 0$, define $b_n =\eta^{}_n$ and
$a^{}_{n,\ell n} = 0$ for all $\ell\geq 1$.
If $\eta^{}_n < 0$, set $a^{}_{n,2n}=-\eta^{}_n$, $a^{}_{n,\ell n}=0$ 
for all $\ell \not= 2$ and $b_n = 0$. 
For $n$ even and $\eta^{}_n-a^{}_{n/2,n}\geq 0$, set
$b_n =\eta^{}_n - a^{}_{n/2,n}$, $a^{}_{d,n} = 0$ for all $d|n$;
for $\eta^{}_n-a^{}_{n/2,n}< 0$ set
$a^{}_{n,2n} = -(\eta_n - a^{}_{n/2,n})$, $a^{}_{n,\ell n}= 0$
for all $\ell \not= 2$ and $b_n = 0$.
Thus \eqref{nthCoeff} yields
$b_n + a^{}_{n/2,n}-a^{}_{n,2n}$ as the $n$-th coefficient,
which, in each of the cases considered above, coincides
with $\eta^{}_n$.
Hence, $H(x)$ can be written as in Proposition~\ref{powerSeriesCrit}.
\end{proof}
\begin{rem}
Analogous to the proof of the Basic Lemma \cite{purisDiss},
Proposition~\ref{powerSeriesCrit} and Theorem~\ref{BLanalogue}
together give a construction of a realizing permutation system for a
given integer sequence.
\end{rem}
\begin{rem}
A combinatorially easier construction can be obtained 
by defining an infinite preimage-cycle for each cycle of the factor
system $(Y,S)$ and, on the other hand, providing $(Y,S)$ with some
$1$-orbit to which all orbits of $(X,T)$ 
can be sent. The drawback of such a construction is
that, due to the infinite preimages of the periodic orbits in the
factor system,
there is no direct way of obtaining a result about the existence of
a realization by compact dynamical systems as given in the Basic
Lemma. 
In this context, it makes sense to introduce a notion like
`factor surjective on the set of periodic points' referring to the
restriction of the factor map to the sets of periodic points being
surjective.
In other words, this defines a
subclass of dynamical systems with factors in which every periodic
orbit of the factor has a periodic orbit in its preimage.
Since the permutation systems in Proposition~\ref{powerSeriesCrit}
are of that type, it follows  that requiring dynamical systems to be
factor surjective on the set of periodic points
is not restrictive with regard to relative realizability.
\end{rem}
\begin{corollary}\label{ERsubsetRR}
    Every exactly realizable sequence is relatively realizable: 
    $\mathcal{ER} \subset \relReal$.
\end{corollary}
\begin{rem}
Obviously, the criterion for exact realizability is
subsumed by the one for relative realizability, but it is
also easy to construct a relative realization for a sequence from
$\mathcal{ER}$: 
for $f$ realized by $(X,T)$, the sequence $2f$ is realized by the 
induced mapping on the topological sum $X+X$ whose factor is the
original map, resulting in a relative realization of the 
sequence $f$.
A further construction is given by the exact realization of 
a sequence $f+u$, where $u_n=1$ for all $n$
and the trivial system $(\{0\},id)$, being a factor of
any dynamical system.
\end{rem}
The following result shows that $\relReal$ shares many properties
with $\mathcal{ER}$ (cf.\ \cite[Section~2]{purisDiss}).
The proofs are either based on the integrality condition of the orbit
counts sequences or
on the construction of realizing dynamical systems. 
They are very similar to those for $\mathcal{ER}$ \cite{purisDiss} and can be
found in \cite{NN}.
\begin{theorem}
    The set $\relReal$ of relatively realizable
    sequences satisfies the following properties:\hfill
\begin{enumerate}
    \item There are no zero divisors in $\relReal$:  $fg = 0$ for
	  $f,g\in\relReal$ implies $f=0$ or
	  $g=0$.\label{zeroDiv}
    \item $\relReal$ contains the constant sequences over $\ZZ$.\label{constant}
    \item $\relReal$ is closed under addition, multiplication
	  und multiplication with elements $z\in\ZZ$.\label{closed}
    \item The constant sequence $u = (1)_{n\geq 1}$ is the only
	  completely multiplicative sequence in $\relReal$.\label{multiplicative}
    \hfill$\Box$
\end{enumerate}
\end{theorem}
\begin{example}\label{TMPD}
The \textit{Thue-Morse}- and \textit{Period Doubling}-chains
are examples of inflation dynamical systems that 
can be treated symbolically.
For background information, consider the articles by Queffelec
\cite{quasicrystalsAF}, 
Allouche and Mend\`es France \cite{quasicrystalsQ} and the references given there.
Define, over the finite alphabets $\{a,b\}$ and $\{A,B\}$, respectively,
the substitution rules
\begin{equation*}
    \mathrm{TM}: a\rightarrow ab,\ b\rightarrow ba\ \text{and}\
    \mathrm{PD}: A\rightarrow AB,\ B\rightarrow AA.
\end{equation*}
The squares $\mathrm{TM}^2$ and $\mathrm{PD}^2$ of each of the two
mappings, iterated on the respective one-letter seeds, produce the
fixed points
$\{.. a|a.. , ..b|b..\}$, $\{.. a|b.., ..b|a..\}$ 
and
$\{..A|A..,..B|B..\}$, $\{..B..,..A..\}$, respectively, i.e., 
$2$-cycles of the original maps.
The inflations $\mathrm{TM}$ and $\mathrm{PD}$ define maps on the 
corresponding $\mathrm{LI}$ classes 
which coincides with the hulls that are obtained as the orbit
closures
under the continuous translation action of $\RR$.

The factor map $\phi$ is given by the block map
\begin{equation*}
    aa\mapsto B,\ bb\mapsto B,\ ab\mapsto A,\ ba\mapsto A.
\end{equation*}
Since the sequences $.. A|A ..$ and $.. B|A..$
correspond to the two sequences $\{.. a|b.., .. b|a.. \}$
and $\{.. a|a.., .. b|b..\}$, respectively, the map $\phi$ is
strictly $2:1$.  (The vertical line indicates the reference point of
the bi-infinite chain.) Following the method described by Anderson
and Putnam
\cite{putnam}, the dynamical zeta function of both systems can be
calculated as \cite{franz}
\begin{equation*}
    \zeta (z) = \frac{1-z}{(1+z)(1-2z)},
\end{equation*}
giving rise to the fixed point counts sequence 
$a_m = 2^m + (-1)^m -1$ (A099430 in the OEIS \cite{OEIS}).
Thus, these two dynamical systems share the same fixed point counts
and relatively realize the sequence
$(0,0,0,\ldots)$ (A000004).
In fact, this and the property of being $2:1$ completely determine
the combinatorics of the map $\phi$.
Let 
\begin{align*}
    \alpha(n) &\coloneqq\ \text{\# $n$-orbits of PD whose preimages are 
				two $n$-orbits of TM}\ \text{and}\\
    \beta(n) &\coloneqq\ \text{\# $n$-orbits of PD whose preimages are
				2$n$-orbits of TM}.
\end{align*} 
Clearly, $c_n = \alpha(n) + \beta(n)$. 
Furthermore, $\beta(n) =\frac{c_n}{2}$ if $n$ is odd and
$\beta(n) = (c_n + \beta(n/2))/2$ if $n$ is even.
Since $a_n$ and $c_n$ are related via \eqref{E:fixOrb}, the
calculation of $\beta(n)$ yields
\begin{equation*}
    \beta(n) = \frac{1}{2n} \sum_{d|n,\ d\ \text{odd}} \mu(d)\cdot 2^{n/d} -
		\delta_{n,1},
\end{equation*}
where $\delta_{n,1} = 1$ for $n=1$ and $0$ otherwise.
Except for $\beta(1) = 0$, this is A000048.
\end{example}
\begin{example}\label{fiboTorus}
The \textit{torus parametrization} of substitution tilings 
\cite{torusParam} yields a large class of dynamical systems with
a torus automorphism as a topological factor.
The one-dimensional Fibonacci chain, obtained by the standard
projection method,
gives rise to the relative
fixed point counts sequence $h_n = (-1)^n$ (arising from the 
fixed point counts A001610 and A001350), which corresponds to
an orbit counts sequence of $\orb(h)=(-1,1,0,0,\ldots)$.
More complicated examples are provided by higher dimensional
projections. 
The Penrose tiling and its torus parametrization, for
instance, lead to the relative
fixed point and orbit counts sequences 
\begin{equation*}
    (-1,9,-16,29,-51,84,-141,\ldots)\ \text{and}\
    (-1,5,-5,5,-10,15,-20,\ldots)
\end{equation*}
which can be calculated from the corresponding dynamical zeta
functions stated explicitly by Baake and Grimm \cite{penroseZeta}.
The more complicated difference sequences in the last case reflect
the phenomenon of a large (though always finite) number of singular
tilings being sent to the same torus parameter ocurring in higher 
dimensional systems, 
whereas the first case illustrates that the one-dimensional Fibonacci  
torus parametrisation is `nearly one-to-one'.
\end{example}
Another large group of  dynamical systems with non-trivial factors
is provided by $S$-integer dynamical systems, cf.\ \cite{Sintegers},
\cite{combinatorialRank} and
references given therein.
Let $Q, P$ be subsets of the set of all primes with $Q \subset P$.
Then the map $S$, dual to $x \mapsto 2x$ on the ring of $Q$-integers,
is a factor of $T$, the map dual to $x \mapsto 2x$ on the ring of
$P$-integers via the dual of the inclusion map from the ring of 
$P$-integers to the ring of $Q$-integers.
\begin{example}\label{SIntegers}
Choose $P$ to be the set of all rational primes, $Q=\varnothing$ and
$\alpha:\QQ\to\QQ$, $x\mapsto 2x$. The
resulting commutative diagramme is obtained by setting $X=\hat{\QQ}$
and $Y=\hat{\ZZ}$ in \eqref{CD}.
The dual map $\hat{\alpha}:\hat{\QQ}\to \hat{\QQ}$ is characterized by
$\chi \mapsto \chi\circ\alpha$ for all $\chi\in\hat{\QQ}$.
Thus, the fixed point equation $\hat{\alpha}^k(\chi) = \chi$ is
equivalent with $\chi(2^k x) = \chi(x)$ for all $x\in\QQ$ or, by the
properties of characters, $\chi((2^k - 1)x) = 1$ for all $x\in\QQ$.
An appropriate choice of $x$ shows that only the trivial
character satisfies this condition,
yielding a fixed point counts sequence of
$(1,1,1,1,\ldots)$ (A000012).
On $Y$, each $k\in\ZZ$ gives rise to the element $e^{2\pi i
k}\in\hat{\ZZ}$, such that the $(2^n-1)$-th roots of unity 
constitute the $2^n - 1$ fixed points of $S^n$ (A000225).

For $P$ chosen as above and $Q =
P\setminus\{3\}$ we obtain $X=\hat{\QQ}$ and $Y=\widehat{\ZZ_{(3)}}$,
where $\ZZ_{(3)} = \ZZ[\frac{1}{p} : p\ \text{prime},\ p\not=
3]$. According to \cite[Example~4.1]{combinatorialRank}, the  
$n$-periodic points of $S$ are given by $\abs{2^n-1}^{-1}_3$, i.e.,
$(1,3,1,3,1,9,1,3,1,3,1,9,1,3,1,\ldots)$, yielding a relative fixed
point sequence of
$( 0, -2, 0, -2, 0, -8, 0, -2, 0, -2, 0, -8, 0, -2, 0,\ldots)$.
\end{example}
\subsection*{Acknowledgements}
It is my pleasure to thank Michael Baake, Franz G\"ahler, Christian Huck 
and Tom Ward for discussions and helpful suggestions.
This work was supported by DFG, within the CRC 701.

\bigskip
This paper is concerned with the integer sequences
\texttt{A000004}, \texttt{A000007}, \texttt{A000012},
\texttt{A000048}, \texttt{A001610}, \texttt{A001350},
\texttt{A060280} and \texttt{A099430} 
from \cite{OEIS}.
\end{document}